\documentclass[french,english]{amsart}
\PassOptionsToPackage{obeyFinal}{todonotes}
\usepackage[T1]{fontenc}
\usepackage[latin9]{inputenc}
\usepackage{xcolor}
\usepackage{babel}
\makeatletter
\addto\extrasfrench{%
   \providecommand{\fg}{\ifdim\lastskip>\z@\unskip\fi~\frqq}%
}

\makeatother
\usepackage{mathtools}
\usepackage{todonotes}
\usepackage{amsbsy}
\usepackage{amstext}
\usepackage{amsthm}
\usepackage{amssymb}
\usepackage{stmaryrd}
\usepackage{graphicx}
\usepackage[unicode=true,pdfusetitle,
 bookmarks=true,bookmarksnumbered=false,bookmarksopen=false,
 breaklinks=false,pdfborder={0 0 1},backref=false,colorlinks=false]
 {hyperref}

\makeatletter

\providecommand{\tabularnewline}{\\}

\numberwithin{equation}{section}
\numberwithin{figure}{section}
\theoremstyle{plain}
\newtheorem{thm}{\protect\theoremname}[section]
\newtheorem*{mainthm}{Main Theorem}
\theoremstyle{definition}
\newtheorem{defn}[thm]{\protect\definitionname}
\theoremstyle{remark}
\newtheorem{notation}[thm]{\protect\notationname}
\theoremstyle{plain}
\newtheorem{prop}[thm]{\protect\propositionname}
\theoremstyle{remark}
\newtheorem{rem}[thm]{\protect\remarkname}
\theoremstyle{plain}
\newtheorem{lem}[thm]{\protect\lemmaname}
\theoremstyle{remark}
\newtheorem{claim}[thm]{\protect\claimname}
\theoremstyle{plain}
\newtheorem{cor}[thm]{\protect\corollaryname}
\theoremstyle{definition}

\usepackage{xcolor}

\DeclareMathOperator*{\Acc}{Acc}

\makeatother

\addto\captionsenglish{\renewcommand{\claimname}{Claim}}
\addto\captionsenglish{\renewcommand{\corollaryname}{Corollary}}
\addto\captionsenglish{\renewcommand{\definitionname}{Definition}}
\addto\captionsenglish{\renewcommand{\examplename}{Example}}
\addto\captionsenglish{\renewcommand{\lemmaname}{Lemma}}
\addto\captionsenglish{\renewcommand{\notationname}{Notation}}
\addto\captionsenglish{\renewcommand{\propositionname}{Proposition}}
\addto\captionsenglish{\renewcommand{\remarkname}{Remark}}
\addto\captionsenglish{\renewcommand{\theoremname}{Theorem}}
\addto\captionsfrench{\renewcommand{\claimname}{Affirmation}}
\addto\captionsfrench{\renewcommand{\corollaryname}{Corollaire}}
\addto\captionsfrench{\renewcommand{\definitionname}{Définition}}
\addto\captionsfrench{\renewcommand{\examplename}{Exemple}}
\addto\captionsfrench{\renewcommand{\lemmaname}{Lemme}}
\addto\captionsfrench{\renewcommand{\notationname}{Notation}}
\addto\captionsfrench{\renewcommand{\propositionname}{Proposition}}
\addto\captionsfrench{\renewcommand{\remarkname}{Remarque}}
\addto\captionsfrench{\renewcommand{\theoremname}{Théorème}}
\providecommand{\claimname}{Claim}
\providecommand{\corollaryname}{Corollary}
\providecommand{\definitionname}{Definition}
\providecommand{\examplename}{Example}
\providecommand{\lemmaname}{Lemma}
\providecommand{\notationname}{Notation}
\providecommand{\propositionname}{Proposition}
\providecommand{\remarkname}{Remark}
\providecommand{\theoremname}{Theorem}

\begin{document}
\global\long\def\R{\mathbb{R}}%

\global\long\def\C{\mathbb{C}}%

\global\long\def\K{\mathbb{K}}%

\global\long\def\Q{\mathbb{Q}}%

\global\long\def\Z{\mathbb{Z}}%

\global\long\def\N{\mathbb{N}}%

\global\long\def\M{\mathcal{M}}%

\global\long\def\U{\mathbb{U}}%

\global\long\def\e{\text{e}}%

\global\long\def\A{\mathcal{A}}%

\global\long\def\B{\mathcal{B}}%

\global\long\def\X{\mathcal{X}}%

\global\long\def\norm#1#2{\left\Vert #1\right\Vert _{#2}}%

\global\long\def\diff{\,d}%

\global\long\def\tendsto#1#2{\underset{#1\to#2}{\longrightarrow}}%

\global\long\def\transpose#1{\prescript{t}{}{#1}}%

\global\long\def\weakto#1#2{\underset{#1\to#2}{\rightharpoonup}}%

\global\long\def\Ml{\M_{0}^{l}}%

\global\long\def\Mlg{\M_{0}^{l,G}}%

\global\long\def\Mli{\M_{0}^{l,I}}%

\global\long\def\Ms{\M_{0}^{s}}%

\global\long\def\Mas{\M_{\text{asym}}}%

\global\long\def\Msrb{\M_{\text{SRB}}}%

\global\long\def\uni{\mathrm{Uni}}%

\global\long\def\erg{\mathrm{Erg}}%

\global\long\def\maj{\mathrm{maj}}%

\global\long\def\floor#1{\left\lfloor #1\right\rfloor }%

\global\long\def\chapeau#1{\text{\ensuremath{\widehat{\delta_{#1}}}}}%

\title[Flexibility versus genericity of phase diagrams]{Flexibility versus genericity of phase diagrams of perturbed continuous maps on the Cantor set}

\author{Hugo MARSAN}
\address{Institut de Math\'ematiques de Toulouse, UMR5219; Universit\'e de Toulouse; CNRS; UPS, F-31062 Toulouse Cedex 9, France.}

\email{hugo.marsan@ens-paris-saclay.fr}

\author{Mathieu SABLIK}
\address{Institut de Math\'ematiques de Toulouse, UMR5219; Universit\'e de Toulouse; CNRS; UPS, F-31062 Toulouse Cedex 9, France.}
\thanks{M.~Sablik acknowledges the support of the ANR ''Difference'' project (ANR-20-CE40-0002).}

\begin{abstract}
Consider the dynamical system constitued by a continuous function $F:\A^\N\longrightarrow\A^\N$ where $\A$ is a finite alphabet. The perturbed counterpart, denoted by $F_\epsilon$, is obtained after each iteration of $F$ by modifying each cell independently with probability $\epsilon\in[0,1]$ and choosing the new value uniformly. We characterize the possible sets of $\epsilon \in[0,1]$ such that $F_\epsilon$ has a unique measure. These sets are exactly the $G_\delta$ sets (countable intersection of open sets) of $[0,1]$ which contain 1. However, we show that 
generically this set is $]0,1]$.
\end{abstract}

\maketitle

A cellular automaton is a continuous function on the Cantor set, either $\A^\N$ or $\A^{\Z^d}$ where $\A$ is a finite alphabet, which commutes with the shift. It is natural to study their perturbed counterpart called Probabilistic Cellular Automata. After each iteration of a given cellular automaton, we modify each cell independently with probability $\epsilon$ and the new value is uniformly chosen over the alphabet $\A$. 

The most important question regarding the asymptotic behaviour of a probabilistic cellular automaton concerns its ergodicity. A probabilistic cellular automaton is said to be ergodic if  its action on probability measures has a unique fixed point that attracts all the other measures. This means that it asymptotically ''forgets'' its initial condition since the distribution of the initial configuration always converges to the same distribution.  It is not known whether a perturbed cellular automaton with a positive rate and only one invariant measure can be not ergodic. However, there is an example of a probabilistic cellular automaton with a non-positive rate which has only one invariant measure but which is not ergodic~\cite{Chassaing-Mairesse-2011}.

Using a percolation argument, a perturbed cellular automaton is ergodic for sufficiently large noise~\cite{MST19}. In other words, the cellular automaton cannot prevent the initial configuration from being forgotten. Moreover large classes of perturbed cellular automata are ergodic~\cite{Vasilyev-78,Gray-1982,MST19,Gacs-Torma-2022}. Constructing a cellular automaton robust to noise in the sense that its trajectories remain distinguishable under the influence of noise is a notoriously difficult problem. The first examples of robust CA are given by A.~Toom~\cite{Toom80} in two-dimensional space. In dimension one, P. Gacs proposes  a highly complex example~\cite{Gacs01}. In all these constructions, the perturbed cellular automata become non ergodic for sufficiently small noise. This means that there is at least one phase transition when the proportion of noise, denoted by $\epsilon$, is considered as a parameter. Recently, a perturbed cellular automaton with at least two phase transitions is shown to exist~\cite{Marsan-Sablik-Torma-2025}, and a natural question is to determine which types of phase transition are possible.

Working with cellular automata is quite difficult, so we transpose this question to continuous functions on $\A^\N$ which do not necessary commute with the shift. Given a continuous function $F: \A^\N\to\A^\N$, we define the perturbed version $F_\epsilon$ where, after the iteration of $F$, each cell is independently modified with probability $\epsilon$  by choosing a letter uniformly at random from $\A$. The aim is to identify the possible phase diagram, i.e. the level of noise at which there is a single invariant measure, as opposed to several. More formally, we want to characterize the possible sets that can be obtained as
$$\uni\left(F\right)=\left\{\epsilon\in[0,1]:F_\epsilon \textrm{ has a unique invariant measure}\right\}.$$

This can be seen as part of the flexibility program launched by Katok, who proposes for a fixed class of dynamical systems to understand the values that a given observable can take~\cite{Erchenko-Katock-2019,Bochi-Katok-2022}. In the thermodynamic formalism on Cantor space, an important question is to characterize the temperature where there is no uniqueness of the measure that maximises pressure. In this setting the flexibility of the phase diagram is explored in~\cite{Kucherenko-Quas-2022}.

Regarding our problem on $\uni\left(F\right)$, the obstruction that arises is that it is a $G_\delta$ containing 1; we detail this in the Section~\ref{Sec.Contraints}. Conversely, in the case of perturbations of continuous Cantor functions which do not necessarily commute with the shift map, any $G_\delta$ of $[0, 1]$ containing 1 can be obtained as $\uni\left(F\right)$ for some continuous function $F:\A^\N\longrightarrow\A^\N$ (Theorem~\ref{thm:Ergodicit=0000E9GDelta}). This new context, where different cells may have neighborhoods of different cardinality, offers greater flexibility, making it possible to prove realisation results much more simply using constructions involving applications of the majority function on well-chosen sets. 
 
The construction uses the majority function introduced in  Section~\ref{Sec.majority} which has two invariant measures for each $\epsilon\in [0,\frac{1}{3}]$ and is simpler than a non-ergodic cellular automaton. Section~\ref{Sec.AuxilliaryLayer}  introduces a new layer that is merely the projection onto $0^\omega$. This is done to retrieve the parameter $\epsilon$, which appears with the noise, by using larger and larger samples to approximate the parameter. This allows to realize  any open set of $\left[0,\frac{1}{3}\right]$ as $\uni\left(F\right)$. To extend the domain in which the function cannot be ergodic, we generalize the majority function in Section~\ref{Sec.majorityVote}. This enables us to realise any open set containing 1 as $\uni\left(F\right)$ in Section~\ref{Sec.RealisationOpen}.  Finally, in Section~\ref{Sec.ConstructionGdelta}, we superimpose the previous construction for a countable number of open sets to obtain the realization result. 
 
Having characterized the possible phase diagrams, we want to examine whether some behaviors emerge. By modifying the previous construction, we can demonstrate that the set of functions realising any $G_\delta$ containing $1$ but not $0$ is dense. Thus the different phase diagrams are widely distributed across the set of functions. However, we can say that a behavior has emerged if it is verified for a generic set, i.e. containing a dense $G_\delta$. In section~\ref{sec.Generic}, we show that a generic dynamical system has a single invariant measure when it is perturbed by a positive noise. It is standard practice to investigate which properties of a class of dynamical systems are generic. For example, see~\cite{Hochman-2008} for different dynamical properties of a class of dynamical system on Cantor set.

The main results of this article can be summarised in the following theorem:
\begin{mainthm}
Let $G$ be a $G_\delta$ of $[0,1]$ which contains 1. There exists $F\in\mathcal{C}(\A^\N)$ such that $\uni\left(F\right)=G$. Moreover, the following set is dense in $(\mathcal{C}(\A^\N),d_{\infty})$:
$$\left\{F\in\mathcal{C}(\A^\N):\uni\left(F\right)=G\setminus\{0\}\right\}.$$

The following set is generic in $(\mathcal{C}(\A^\N),d_{\infty})$:
$$\left\{F\in\mathcal{C}(\A^\N):\uni\left(F\right)=]0,1]\right\}.$$
\end{mainthm}

\section{Definitions}\label{Sec.definition}

\subsection{Set of configurations}
In this article, $\A$ denotes a finite alphabet of symbols. If $I_{n}$ denotes the subset $\left\{ 0,\dots,n\right\} $, then the product topology on $\A^{\N}$ is induced by the distance
\[
d\left(x,y\right)\coloneqq2^{-\min\left\{ n\in\N\,\mid\,x_{I_{n}}\neq y_{I_{n}}\right\} },
\]
which makes $\A^{\N}$ compact. For a finite word $\omega\in\A^{r}$ and $\U=\left\{ u_{1},\dots,u_{r}\right\} \subset\Z^{d}$, the cylinder $\left[\omega\right]_{\U}$ denotes the set of configurations that agree with $\omega$ on $\U$:
\[
\left[\omega\right]_{\U}\coloneqq\left\{ x\in\X\mid x_{\U}=\omega\right\} .
\]
The set of all cylinders is a base of open subsets for the product topology on $\A^{\N}$, and thus induce the Borel $\sigma$-algebra $\mathfrak{B}$.

\subsection{Continuous functions of $\protect\A^{\protect\N}$ and CA}

By compactness of $\A^\N$, continuous functions of $\A^\N$ can be considered as a function on $\A^\N$ where the updated value in each cell depends on a finite neighborhood, whose size   depends on the cell. 
\begin{notation}
Given a continuous function $F$, denote by $\left(\mathcal{N}_{i}\right)_{i\in\N}$ its local neighborhoods and $\left(f^{i}\right)_{i\in\N}$ its local rules, such that for each $i\in\N$, $\mathcal{N}_{i}$ is a finite subset of $\N$ and $f^{i}:\A^{\mathcal{N}_{i}}\to\A$ that verify for all $x\in\A^{\N}$, 
\[
\left(F\left(x\right)\right)_{i}=f^{i}\left(x_{\mathcal{N}_{i}}\right).
\]
\end{notation}

If the update functions are the same, the continuous function is a Cellular Automaton, equivalently is a continuous function which commutes with the shift~\cite{Morse-Hedlund-1938}. 

\begin{defn}
A function $F:\A^\N\to\A^\N$ is a Cellular Automaton if there exists a finite $\mathcal{N}\subset\N$ and $f:\A^{\mathcal{N}}\to\A$ such that for all $i\in\N$, $\mathcal{N}_{i}=i+\mathcal{N}=\left\{ i+k\mid k\in\mathcal{N}\right\} $ and $f^{i}\equiv f$. They are exactly the continuous function that commutes with the shift $\sigma:\left(x_{i}\right)_{i\in\N}\longmapsto\left(x_{i+1}\right)_{i\in\N}$.
\end{defn}

The set of continuous function is denoted $ \mathcal{C}\left(\A^{\N}\right)$. This space becomes a compact metric space when it is endowed with the distance $$d_\infty(F,F')=\sup_{x\in\A^\N}d(F(x),F'(x)).$$

\subsection{Probability measures and perturbations}

A sequence $\left(\mu_{n}\right)_{n\in\N}$ of probability measures of $\A^{\N}$ is said to weakly converge towards $\mu$ (denoted by $\mu_{n}\weakto n{\infty}\mu$) if for all cylinders $\left[\omega\right]_{\U}$, $\mu_{n}\left(\left[\omega\right]_{\U}\right)\tendsto n{\infty}\mu\left(\left[\omega\right]_{\U}\right)$. The set $\M\left(\A^{\N}\right)$ is compact for the induced weak converge topology.
\begin{defn}
\label{def:MesuresProba}Fix $\U\subset\N^{d}$ finite and $\omega\in\A^{\U}$.
\begin{itemize}
\item $\lambda$ denotes the uniform Bernoulli measure, defined by 
\[
\lambda\left(\left[\omega\right]_{\U}\right)\coloneqq\frac{1}{\left|\A^{\U}\right|}.
\]
\item For $\alpha\in\M\left(\A\right)$ ($\alpha=\left(\alpha_{b}\right)_{b\in\A}$ such that $\sum\alpha_{b}=1$ and $\alpha_{b}\geq0$), $\lambda_{\alpha}$ denotes the Bernoulli measure with parameter $a$, defined by
\[
\lambda_{\alpha}\left(\left[\omega\right]_{\U}\right)\coloneqq\prod_{i\in\U}\alpha_{\omega_{i}}.
\]
\end{itemize}
\end{defn}

\begin{notation}
In the case $\A=\left\{ 0,1\right\} $, we also use the notation $\lambda_{x}$ with $x\in\left[0,1\right]$ for $\lambda_{\left(1-x,x\right)}$. The uniform measure on $\left\{ 0,1\right\} ^{\N}$ can then be written $\lambda$, $\lambda_{\frac{1}{2}}$ or $\lambda_{\left(\frac{1}{2},\frac{1}{2}\right)}$.
\end{notation}

\begin{defn}
Let $F\in \mathcal{C}\left(\A^{\N}\right)$ with local neighborhoods $\left(\mathcal{N}_{i}\right)_{i\in\N}$ and local rules $\left(f^{i}\right)_{i\in\N}$. For $\epsilon\geq0$, define
\begin{align*}
f_{\epsilon}^{i}:\A^{\mathcal{N}_{i}}\times\A & \to\left[0,1\right]\\
\left(a,b\right) & \mapsto\left(1-\epsilon\right)\mathbf{1}_{b=f^{i}\left(a\right)}+\frac{\epsilon}{\left|\A\right|}
\end{align*}

The \emph{computer perturbation} of $F$ of size $\epsilon$ is the probability kernel $F_{\epsilon}:\A^{\N}\times\mathfrak{B}\to\left[0,1\right]$ with local rules $\left(f_{\epsilon}^{i}\right)$, i.e. verifying for all $x\in\A^{\N}$, $\U\subset\N$ and $\omega\in\A^{\U}$, 
\[
F_{\epsilon}\left(x,\left[\omega\right]_{\U}\right)=\prod_{i\in\U}f_{\epsilon}^{i}\left(x_{\mathcal{N}_{i}},\omega_{i}\right).
\]
\end{defn}

\begin{rem}
The noise is called computer perturbation since each bit is independently perturbed with a certain probability, see for example~\cite{AsaCol05} which propose a model of perturbed Turing machine.  By contrast, physical noise changes a bit with a probability that increases as the bit is further away from the origin.
\end{rem}

\begin{prop}
For $F\in \mathcal{C}\left(\A^{\N}\right)$ and $\epsilon\in\left[0,1\right]$, the action of $F_{\epsilon}$ on $\M\left(\A^{\N}\right)$ is defined by 
\[
F_{\epsilon}\mu\left(A\right)\coloneqq\int_{x\in\A^{\N}}F_{\epsilon}\left(x,A\right)\diff\mu\left(x\right)
\]
for all $A\in\mathfrak{B}$, is continuous.
\end{prop}

\begin{rem}
At $\epsilon=0$, the action is defined by $F\mu\left(A\right)=\mu\left(F^{-1}\left(A\right)\right)$.
\end{rem}

By compactness of $\M(\A^\N)$, the function $F_\epsilon$ admits invariant measures.  Denote this set by $\M^F_\epsilon$, or simply $\M_\epsilon$ if there is no ambiguity. 

\subsection{Ergodic theory}

A probability measure $\mu$ is said to be $F_{\epsilon}$-invariant if $F_{\epsilon}\mu=\mu$. By continuity of the action of $F_{\epsilon}$ and compactness of $\M\left(\A^{\N}\right)$, the set of invariant measures of $F_{\epsilon}$, denoted by $\M_{\epsilon}$, is a compact and non-empty set. The following definition takes its vocabulary from the theory of probabilistic cellular automata.
\begin{defn}
For $F\in \mathcal{C}\left(\A^{\N}\right)$ and $\epsilon\geq0$, $F_{\epsilon}$ is said to be ergodic if its only admits a unique invariant measure, which is also attractive: there exists a probability measure $\pi\in\M\left(\A^{\N}\right)$ such that $\M_{\epsilon}=\left\{ \pi\right\} $ and for all $\mu\in\M\left(\A^{\N}\right)$, 
\[
F_{\epsilon}^{n}\mu\weakto n{\infty}\pi.
\]
\end{defn}

Finding a cellular automaton that is not ergodic in the face of this kind of perturbation is a challenging problem. A. Toom~\cite{Toom80} gives an exemple in two dimension and P. Gacs proposes  a highly complex example~\cite{Gacs01} in dimension one.

\begin{defn}
Let $F$ be a continuous maps of $\A^{\N}$ and $F_{\epsilon}$ its computer perturbation of size $\epsilon$. The two sets $\uni\left(F\right)$ and $\erg\left(F\right)$ are defined by 
\begin{align*}
\uni\left(F\right) & =\left\{ \epsilon\in\left[0,1\right]\mid F_{\epsilon}\text{ admits a unique invariant measure}\right\} \\
\erg\left(F\right) & =\left\{ \epsilon\in\left[0,1\right]\mid F_{\epsilon}\text{ is ergodic}\right\} .
\end{align*}
\end{defn}

By definition, $\erg\left(F\right)\subset\uni\left(F\right)$. As the number of known examples of non-ergodic 1-dimensional probabilistic cellular automata with positive rate is so low, the question of whether there is a cellular automaton $F$ such that $\erg\left(F\right)\neq\uni\left(F\right)$ is still open (see for example \cite{CM11} for an example without positive rates).

\section{Constraints for ergodicity}\label{Sec.Contraints}

In this section we explore the obstruction for $\uni\left(F\right)$ when $F$ is a continuous function.

\subsection{Continuity lemma}

To exhibit topological contraints we need a continuity lemma according the parameter $\epsilon$. In particular when $\epsilon$ goes to $\epsilon_0$, we obtain that the accumulation points of thes sets of invariant measure for the level of noise $\epsilon$, denoted $\M_{\epsilon}$,  is included in $\M_{\epsilon_0}$. This Lemma is a generalisation for the continuous fonction of the same lemma for cellular automata obtained in~\cite{Marsan-Sablik-2024}

\begin{lem}
\label{lem:ContinuityLemma} Let $(F^i)_{i\in\N}$ be a sequence of continuous functions of $\A^{\N}$ which converges to $F$ and  $(\epsilon_i)_{i\in\N}$ be a sequence of elements of $[0,1]$ which converges to~$\epsilon$. Let $(\pi_{i})_{i\in\N}$ be a family of measures such that $\pi_{i}\weakto{i}{\infty}\pi$. Then $$F^i_{\epsilon_i}\pi_{i}\weakto{i}{\infty}F_\epsilon\pi.$$

In particular, $\Acc_{\epsilon\to\epsilon_{0}}\left(\M_{\epsilon}\right)\subset\M_{\epsilon_{0}}$.
\end{lem}

\begin{lem}
\label{lem:ContinuityLemma} Let $(F^i)_{i\in\N}$ be a sequence of continuous functions of $\A^{\N}$ which converges to $F$ and  $(\epsilon_i)_{i\in\N}$ be a sequence of elements of $[0,1]$ which converges to~$\epsilon$. Let $(\pi_{i})_{i\in\N}$ be a family of measures such that $\pi_{i}\weakto{i}{\infty}\pi$. Then $$F^i_{\epsilon_i}\pi_{i}\weakto{i}{\infty}F_\epsilon\pi.$$

In particular, $\Acc_{\epsilon\to\epsilon_{0}}\left(\M_{\epsilon}\right)\subset\M_{\epsilon_{0}}$.
\end{lem}

\begin{proof}
Let $\U\subset\Z^{d}$ be a finite subset. Let us show that $\norm{F^i_{\epsilon_i}\pi_{i}-F_{\epsilon}\pi}{\U}\tendsto{i}{\infty}0$. By triangular inequality, $\norm{F^i_{\epsilon_i}\pi_{i}-F_{\epsilon}\pi}{\U}\leq\underset{\boldsymbol{2}}{\underbrace{\norm{F^i_{\epsilon_i}\pi_{i}-F_{\epsilon}\pi_{i}}{\U}}}+\underset{\boldsymbol{1}}{\underbrace{\norm{F_{\epsilon}\pi_{i}-F_{\epsilon}\pi}{\U}}}$. 
\begin{enumerate}
\item[\foreignlanguage{french}{$\boldsymbol{1}$}] tends to $0$ by continuity of the action of $F_{\epsilon}$ on $\M\left(\A^{\N}\right)$. 
\item[\foreignlanguage{french}{$\boldsymbol{2}$}] tends to $0$ if for all $u\in\A^{\mathbb{U}}$, $F^i_{\epsilon_i}\pi_{i}\left(\left[u\right]_{\U}\right)-F_{\epsilon}\pi\left(\left[u\right]_{\U}\right)\tendsto{i}{\infty}0$. 

First remark that for $i$ sufficiently large, all the local rules of $F^i$ for the coordinate $j\in\U$ have the same local rule as $F$. Using
\[
\left|F^i_{\epsilon_i}\pi_{i}\left(\left[u\right]_{\U}\right)-F_{\epsilon}\pi\left(\left[u\right]_{\U}\right)\right|\leq\int\left|F^i_{\epsilon_i}\left(x,\left[u\right]_{\U}\right)-F_{\epsilon}\left(x,\left[u\right]_{\U}\right)\right|\diff\pi_{\epsilon},
\]
observe that with $\Delta_{\epsilon_i,\epsilon}\coloneqq\max_{i\in\U}\norm{f_{\epsilon_i}^{j}-f_{\epsilon}^{j}}{\infty}$,
\begin{align*}
F^i_{\epsilon_i}\left(x,\left[u\right]_{\U}\right) & \coloneqq\prod_{j\in\U}f_{\epsilon_i}^{j}\left(x_{\mathcal{N}_{j}},u_{j}\right)\\
 & \leq\prod_{j\in\U}\left(f_{\epsilon}^{j}\left(x_{j+\mathcal{N}},u_{j}\right)+\Delta_{\epsilon_i,\epsilon}\right) & \textrm{}\\
 & \leq\left(\prod_{j\in\U}f_{\epsilon}^{j}\left(x_{j+\mathcal{N}},u_{i}\right)\right)+\sum_{k=1}^{\left|\U\right|}\binom{\left|\U\right|}{k}\Delta_{\epsilon_i,\epsilon}^{k} & \begin{array}{c}
\text{by expansion}\\
\text{and }f_{\epsilon}^{j}\in\left[0,1\right]
\end{array}\\
 & =F_{\epsilon}\left(x,\left[u\right]_{\U}\right)+\left(\Delta_{\epsilon_i,\epsilon}+1\right)^{\left|\U\right|}-1
\end{align*}
and by symmetry one can conclude 
\[
\left|F^i_{\epsilon_i}\left(x,\left[u\right]_{\U}\right)-F_{\epsilon}\left(x,\left[u\right]_{\U}\right)\right|\leq\left(\Delta_{\epsilon_i,\epsilon}+1\right)^{\left|\U\right|}-1.
\]
By definition of $f_{\epsilon}^{j}$, one has $\Delta_{\epsilon_i,\epsilon}\leq\left|\epsilon_i-\epsilon\right|\tendsto{i}{\infty}0$ and finally, 
\[
\left|F^i_{\epsilon_i}\pi_{\epsilon}\left(\left[u\right]_{\U}\right)-F_{\epsilon}\pi\left(\left[u\right]_{\U}\right)\right|\leq\left(\Delta_{\epsilon_i,\epsilon}+1\right)^{\left|\U\right|}-1\tendsto{i}{\infty}0.
\]
\end{enumerate}
Hence the result.
\end{proof}

\subsection{Topological constraint}

We can define constraints on the theoretical possible sets. As we only consider the uniform standard perturbation, we already have a continuity (and computability) hypothesis on the noise.
\begin{prop}
\label{prop:uni_est_un_G_delta}Let $F$ be a continuous function on $\A^{\N}$. Then $\uni\left(F\right)$ is a $G_{\delta}$ subset of $\left[0,1\right]$, i.e. an countable intersection of open subsets of $\left[0,1\right]$. Moreover, $1\in\uni\left(F\right)$.
\end{prop}

\begin{proof}
By definition of the perturbation, one has immediately $1\in\erg\left(F\right)$ as for $\epsilon=1$, for every $\mu$ initial probability measure $F_{\epsilon}\mu=\lambda$. 

For the other part, as $$\epsilon\in\uni\left(F\right)\Longleftrightarrow\text{diam}\left(\M_{\epsilon}\right)=0,$$ the set can be defined as $\uni\left(F\right)=\bigcap_{k\in\N^{*}}A_{k}$ where $$A_{k}=\left\{ \epsilon\in\left[0,1\right]\mid\text{diam}\left(\M_{\epsilon}\right)<\frac{1}{k}\right\}.$$ 

To conclude, we only need to prove that their complementary sets $A_{k}^{c}$ are closed. Let $\epsilon_{n}\tendsto n{\infty}\epsilon\in\left[0,1\right]$ a sequence of $A_{k}^{c}$. By compactness of $\M_{\epsilon_{n}}$, there exists $\mu_{n}$ and $\nu_{n}$ such that $d_{\M}\left(\mu_{n},\nu_{n}\right)\geq\frac{1}{k}$. Using sub-sequences and the compactness of $\M\left(\A^{\N}\right)$, one can suppose that $\mu_{n}$ and $\nu_{n}$ respectively converge to $\mu$ and $\nu$. The distance being continuous, they verify $d_{\M}\left(\mu,\nu\right)\geq\frac{1}{k}$ and by continuity Lemma \ref{lem:ContinuityLemma}, $\mu,\nu\in\M_{\epsilon}$, and thus $\epsilon\in A_{k}^{c}$.
\end{proof}

\section{Realization}

In this section we describe the construction of a continuous function whose computer perturbation is ergodic only when the error rate $\epsilon$ belongs to an arbitrary $G_{\delta}$ set, and admits several invariant measures otherwise. We re-use the main idea of the article~\cite{Marsan-Sablik-Torma-2025}: we first need a continuous function which is robust to perturbations and then couple it with an adversary which can force its ergodicity on a given error rate. In the role of the robust map, we use a kind of majority vote, which was already used as the NEC-majority to find an example of a non-ergodic positive-rate PCA by Toom in~\cite{Toom80}.

We first use a classic majority vote on 3 cells, show its robustness and what we can obtain when coupled with an adversary. We then generalize the model to a majority vote on $2n+1$ cells, to obtain more general ergodicity sets. The last step will be to couple all those majority vote to finally obtain any $G_{\delta}$ sets containing $1$ as an ergodicity set.

\subsection{Majority vote on 3 cells}\label{Sec.majority}
\begin{prop}
\label{prop:maj3}Let $\maj_{3}:\left\{ 0,1\right\} ^{\N}\to\left\{ 0,1\right\} ^{\N}$ be the continuous function defined by:
\[
\maj_{3}\left(x\right)_{i}\coloneqq\maj\left(x_{3i+1},x_{3i+2},x_{3i+3}\right).
\]
Then its computer perturbation of size $0<\epsilon<\frac{1}{3}$, denoted by $\maj_{3,\epsilon}$, admits several invariant measures, including the Bernoulli measures with parameters $\alpha_{\epsilon}$, $\frac{1}{2}$ and $1-\alpha_{\epsilon}$ with $\alpha_{\epsilon}=\frac{1}{2}\left(1-\sqrt{1-\frac{2\epsilon}{1-\epsilon}}\right)$.
\end{prop}

\begin{figure}[h]
\begin{centering}
\includegraphics[width=14cm]{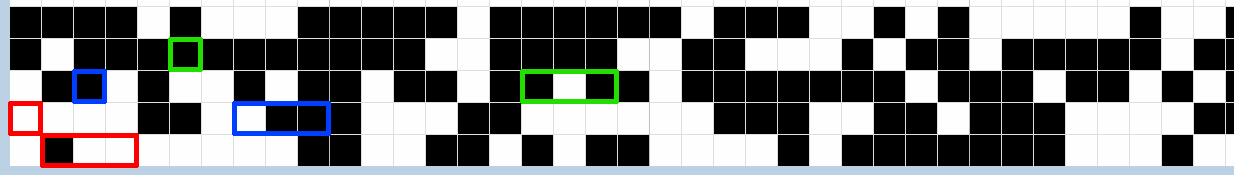}
\par\end{centering}
\caption{Trajectory of $\protect\maj_{3}$, without perturbation. We highlighted certain cells with their neighborhood in the previous iteration in the same color.}
\end{figure}

\begin{proof}
Suppose first that a Bernoulli measure with parameter $\alpha$, denoted by $\lambda_{\alpha}$, is $\maj_{3,\epsilon}$-invariant. Then 
\[
\lambda_{\alpha}\left(\left[1\right]_{0}\right)=\maj_{3,\epsilon}\lambda_{\alpha}\left(\left[1\right]_{0}\right)
\]
which can be rewritten, with $E_{i}$ the event <<an error at $0$ creating the symbol $i$>>:
\begin{align*}
\alpha & =P\left(E_{1}\right)+P\left(\overline{E_{0}\cup E_{1}}\right)\cdot\lambda_{\alpha}\left(\left[011\right]_{\left\llbracket 1,3\right\rrbracket }\cup\left[101\right]_{\left\llbracket 1,3\right\rrbracket }\cup\left[110\right]_{\left\llbracket 1,3\right\rrbracket }\cup\left[111\right]_{\left\llbracket 1,3\right\rrbracket }\right)\\
 & =\frac{\epsilon}{2}+\left(1-\epsilon\right)\left(3\alpha^{2}\left(1-\alpha\right)+\alpha^{3}\right).
\end{align*}
Thus $\alpha$ is a root of the polynomial $P_{\epsilon}=-2(1-\epsilon)X^{3}+3(1-\epsilon)X^{2}-X+\frac{\epsilon}{2}$, so $\alpha\in\left\{ \frac{1}{2},\alpha_{\epsilon},1-\alpha_{\epsilon}\right\} $ with $\alpha_{\epsilon}=\frac{1}{2}\left(1-\sqrt{1-\frac{2\epsilon}{1-\epsilon}}\right)$.

Then, it suffices to show that for all $\alpha$ root of $P_{\epsilon}$, $\lambda_{\alpha}$ is $\maj_{3,\epsilon}$-invariant. By definition, the neighborhood $\mathcal{N}_{i}=\left\{ 3i+1,3i+2,3i+3\right\} $ of each cell is disjoint form the others. Thus the independence between each cell is conserved by $\maj_{3,\epsilon}$: for a finite $\U\subset\N$ and $u\in\left\{ 0,1\right\} ^{\U}$, 
\[
\maj_{3,\epsilon}\lambda_{\alpha}\left(\left[u\right]_{\U}\right)=\prod_{i\in\U}\maj_{3,\epsilon}\lambda_{\alpha}\left(\left[u_{i}\right]_{i}\right).
\]
For all $i\in I$, $\maj_{3,\epsilon}\lambda_{\alpha}\left(\left[1\right]_{i}\right)=\frac{\epsilon}{2}+\left(1-\epsilon\right)\left(3\alpha^{2}\left(1-\alpha\right)+\alpha^{3}\right)=\alpha$ as $\alpha$ is a root of $P_{\epsilon}$. 
One can then obtain $a\in\left\{ 0,1\right\} $, $\maj_{3,\epsilon}\lambda_{\alpha}\left(\left[a\right]_{i}\right)=\lambda_{\alpha}\left(\left[a\right]_{i}\right)$ and finally
\[
\maj_{3,\epsilon}\lambda_{\alpha}\left(\left[u\right]_{\U}\right)=\prod_{i\in\U}\lambda_{\alpha}\left(\left[u_{i}\right]_{i}\right)=\lambda_{\alpha}\left(\left[u\right]_{\U}\right)
\]
which finishes the proof.
\end{proof}

\subsection{Adding an auxiliary layer}\label{Sec.AuxilliaryLayer}

In this section, the alphabet considered is now $\A=\left\{ 0,1\right\} \times\left\{ 0,1\right\} $. The function $\maj_{3}$ is coupled with the constant function equals to $0$. For a fixed error rate $\epsilon$, the density of symbols $1$ on the second layer is exactly $\frac{\epsilon}{2}$, even after a single iteration of the perturbed function. We can use this observation to approximate the error rate and change the behavior of the coupling: when the approximation is in a subset $I_{i}$, a $0$ symbol is projected on the first layer, regardless of the result of $\maj_{3}$.
\begin{defn}
Let $I=\left(I_{i}\right)_{i\in\N}$ be a sequence of subsets of $\left[0,1\right]$, and $E=\left(E_{i}\right)_{i\in\N}$ a sequence of finite subsets of $\N$. For $x=\left(y,z\right)\in\left(\left\{ 0,1\right\} \times\left\{ 0,1\right\} \right)^{\N}$, define $x_{i}\coloneqq\left(y_{i},z_{i}\right)$. $\maj_{3}^{I,E}\in \mathcal{C}\left(\A^{\N}\right)$ denotes the function with local neighborhood $\left\{ 3i,3i+1,3i+2\right\} \cup E_{i}$ and local rule
\[
\maj_{3}^{I,E}\left(x\right)_{i}\coloneqq\begin{cases}
\left(0,0\right) & \text{if }2\cdot\text{mean}\left(z,i\right)\in I_{i}\\
\left(\maj\left(y_{3i+1},y_{3i+2},y_{3i+3}\right),0\right) & \text{otherwise}
\end{cases}
\]
with 
\[
\text{mean}\left(z,i\right)\coloneqq\frac{1}{\left|E_{i}\right|}\sum_{j\in E_{i}}z_{j}
\]
the frequency of symbols $1$ on the second layer of the cells in $E_{i}$.
\end{defn}

\begin{notation}
For $J\subset\N$, define $E_{J}\coloneqq\bigcup_{i\in J}E_{i}$.
\end{notation}

By carefully choosing the families $I$ and $E$, one can obtain a continuous function whose computer perturbation is ergodic only at an arbitrary $\epsilon_{0}<\frac{1}{3}$.
\begin{prop}
\label{prop:maj3IE}Fix $0\leq\epsilon_{0}<\frac{1}{3}$. There exists $E$ and $I$ such that the computer perturbation of size $\epsilon$ of $\maj_{3}^{I,E}$ is ergodic only for $\epsilon=\epsilon_{0}$ when $\epsilon<\frac{1}{3}$, and admits otherwise several invariant measures:
\[
\uni\left(\maj_{3}^{I,E}\right)\cap\left[0,\frac{1}{3}\right[=\erg\left(\maj_{3}^{I,E}\right)\cap\left[0,\frac{1}{3}\right[=\left\{ \epsilon_{0}\right\} .
\]
\end{prop}

\begin{proof}
For $t\geq1$, define $\mathcal{N}^{t}=\mathcal{N}_{0}^{t}\coloneqq\left\llbracket \frac{3^{t}-1}{2},3^{t}+\frac{3^{t}-3}{2}\right\rrbracket $ the local neighborhood of $\left\{ 0\right\} $ by the continuous function $\maj_{3}^{t}$. Observe that $\left(\mathcal{N}^{t}\right)_{t\in\N^{*}}$ is a partition of $\N^{*}$ such that for all $t\geq1$, $3^{t}\in\mathcal{N}^{t}$. In general, if $i\in\mathcal{N}^{k}$, define $\mathcal{N}_{i}^{t}\subset\mathcal{N}^{t+k}$ the local neighborhood of $\left\{ i\right\} $ by $\maj_{3}^{t}$. We can then choose $I$ and $E$ such that:
\begin{itemize}
\item $i\mapsto I_{i}$ is constant on each $\mathcal{N}^{t}$, with value $I_{3^{t}}$.
\item $i\mapsto\left|E_{i}\right|$ is constant on each $\mathcal{N}^{t}$, with value $\left|E_{3^{t}}\right|$. We choose here$\left|E_{3^{t}}\right|=4^{t}$.
\item The $E_{i}$ are pairwise disjoints.
\item For all $i\in\N$ and $t\geq2$ , $\mathcal{N}_{i}^{t-1}$ and $E_{\mathcal{N}_{i}^{t-2}}$ are disjoints (only useful for the non-ergodicity part). It suffices that $\mathcal{N}^{t}\cap E_{\mathcal{N}^{t-1}}=\emptyset$ for all $t\geq1$ to verify this condition.
\end{itemize}
The $E_{i}$ are defined by induction: fix $E_{0}=\left\{ 4\right\} $. If $E_{i}=\left\llbracket a,b\right\rrbracket $ with $i\in\mathcal{N}^{k}$, define $E_{i+1}$ as the following:
\begin{itemize}
\item If $i+1\in\mathcal{N}^{k}$ as well, then $E_{i+1}=\left\llbracket b+1,b+4^{k}\right\rrbracket $.
\item If $i+1\in\mathcal{N}^{k+1}$, then $E_{i+1}=\left\llbracket b^{\prime}+1,b^{\prime}+4^{k+1}\right\rrbracket $ with $b^{\prime}\coloneqq\max\left(b,\max\left(\mathcal{N}^{k+2}\right)\right)$. 
\end{itemize}
By construction, the $E_{i}$ verify the properties previously listed.

\textbf{Ergodicity at $\epsilon=\epsilon_{0}$}: let us first show that $\left(\maj_{3}^{I,E}\right)_{\epsilon_{0}}$ is ergodic. For two initial measures $\mu,\nu\in\M\left(\left(\left\{ 0,1\right\} \times\left\{ 0,1\right\} \right)^{\N}\right)$, define a coupling $\left(X^{t},Y^{t}\right)$ ($X^{t}\sim\mu^{t}\coloneqq\left(\maj_{3}^{I,E}\right)_{\epsilon_{0}}^{t}\mu$ and $Y^{t}\sim\nu^{t}\coloneqq\left(\maj_{3}^{I,E}\right)_{\epsilon_{0}}^{t}\nu$) by applying to each trajectories the same errors: an error creates the symbol $\alpha\in\left\{ 0,1\right\} \times\left\{ 0,1\right\} $ at $i\in\N$ in $X^{t}$ if and only if an error creates the same symbol at $i$ on $Y^{t}$. In particular, as $\maj_{3}^{I,E}$ corresponds to the constant function equals to $0$ on the second coordinate, the second layers of $X^{t}$ and $Y^{t}$ are equals for $t\geq1$: we denote it by $Z^{t}$. One can easily obtain that $Z^{t}\sim Z^{1}\sim\lambda_{\frac{\epsilon_{0}}{2}}$.

Let $t\geq2$. If for all $i\in\mathcal{N}^{t-2}$, $\text{mean}\left(Z^{1},i\right)\in I_{i}$, then $X_{0}^{t}=Y_{0}^{t}$. Indeed, as the same errors are applied to $X$ and $Y$, one has $X_{i}^{2}=Y_{i}^{2}$ for all $i\in\mathcal{N}^{t-2}$. Then, $X_{i}^{k}=Y_{i}^{k}$ for $i\in\mathcal{N}^{t-k}$, as $X$ and $Y$ share the same second layer $Z$ and as for all $j\in\mathcal{N}^{t-k-1}$, $X_{j}^{k+2}$ depends only on the errors at the time $k+1$, on $Z^{k}$ and on $X_{i}^{k}$ with $i\in\mathcal{N}^{t-k}$, then $X_{j}^{k+1}=Y_{j}^{k+1}$: we can conclude by induction on $k$ that $X_{0}^{t}=Y_{0}^{t}$.

By definition, the values $\text{mean}\left(Z^{1},i\right)$ are all independent and the function $i\longmapsto p_{i}\coloneqq P\left(2\cdot\text{mean}\left(Z^{1},i\right)\in I_{i}\right)$ is constant on each $\mathcal{N}^{t}$. Thus one has
\begin{align*}
P\left(X_{0}^{t}=Y_{0}^{t}\right) & \geq P\left(\bigcap_{i\in\mathcal{N}^{t-2}}2\cdot\text{mean}\left(Z^{1},i\right)\in I_{i}\right)=\prod_{i\in\mathcal{N}^{t-2}}p_{i}=\left(p_{3^{t-1}}\right)^{3^{t-2}}.
\end{align*}
Suppose that there exists a sequence $\left(\eta_{t}\right)$ such that for $t\in\N^{*}$ large enough, $I_{3^{t}}\supset\left[\epsilon_{0}\pm\eta_{t}\right]$. Then by Bienaym\'e-Chebychev inequality on $Z^{1}$ with law $\lambda_{\frac{\epsilon_{0}}{2}}$,
\begin{align*}
p_{3^{t}} & \geq1-P\left(\left|\text{mean}\left(Z^{1},3^{t}\right)-\frac{\epsilon_{0}}{2}\right|>\frac{\eta_{t}}{2}\right)\\
 & \geq1-\frac{\frac{\epsilon_{0}}{2}\left(1-\frac{\epsilon_{0}}{2}\right)4}{\left|E_{3^{t}}\right|\eta_{t}^{2}}\\
 & \geq1-\frac{1}{\left|E_{3^{t}}\right|\eta_{t}^{2}}.
\end{align*}
Finally, for $\eta_{t}=\frac{1}{t}$, $I_{3^{t}}=\left[\epsilon_{0}\pm\frac{1}{t}\right]$ and $\left|E_{3^{t}}\right|=4^{t}$, $3^{t}=o\left(\left|E_{3^{t}}\right|\eta_{t}^{2}\right)$ so $\left(p_{3^{t}}\right)^{3^{t-1}}\tendsto t{+\infty}1$ and $P\left(X_{0}^{t}=Y_{0}^{t}\right)\tendsto t{+\infty}1$. More generally, for all $j\in\N$, $i\mapsto I_{i}$ and $i\mapsto\left|E_{i}\right|$ are still constant on each $\mathcal{N}_{j}^{t}$, and similar computations lead to
\[
P\left(X_{j}^{t}=Y_{j}^{t}\right)\geq\left(p_{3^{t-1+k}}\right)^{3^{t-2}}\tendsto t{+\infty}1
\]
so $\left(\maj_{3}^{I,E}\right)_{\epsilon_{0}}$ is ergodic.

\textbf{Non-ergodicity at $\epsilon\neq\epsilon_{0}$}: fix $\epsilon\neq\epsilon_{0}$ (with $\epsilon<\frac{1}{3}$). For initial measure, let us take $\mu$ such that the projection on the first layer is $\lambda_{0}=\chapeau 0$. Denote by $\mu_{\epsilon}^{t}\coloneqq\left(\maj_{3}^{I,E}\right)_{\epsilon}^{t}\mu$ and $\lambda_{\epsilon,\alpha}^{t}\coloneqq\left(\maj_{3}\right)_{\epsilon}^{t}\lambda_{\alpha}$ As $\maj_{3}^{I,E}$ only <<adds>> more 0-symbol on the first layer compared to $\maj_{3}$, we have for $t,j\in\N$ 
\[
\mu_{\epsilon}^{t}\left(\left[\left(1,0\right)\right]_{j}\cup\left[\left(1,1\right)\right]_{j}\right)\leq\lambda_{\epsilon,0}^{t}\left(\left[1\right]_{j}\right).
\]
By the computations of Proposition \ref{prop:maj3}, $\lambda_{\epsilon,0}^{t}=\lambda_{\left(h_{\epsilon}\right)^{t}\left(0\right)}$ with \[h_{\epsilon}\left(\alpha\right)=\frac{\epsilon}{2}+\left(1-\epsilon\right)\left(3\alpha^{2}\left(1-\alpha\right)+\alpha^{3}\right).\] The interval $\left[0,\alpha_{\epsilon}\right[$ with $\alpha_{\epsilon}=\frac{1}{2}\left(1-\sqrt{1-\frac{2\epsilon}{1-\epsilon}}\right)$ is stable by $h_{\epsilon}$, thus $\lambda_{\epsilon,\alpha}^{t}\left(\left[1\right]_{j}\right)\leq\alpha_{\epsilon}$, and then for all $t\in\N$ one has $\mu_{\epsilon}^{t}\left(\left[\left(1,0\right)\right]_{j}\cup\left[\left(1,1\right)\right]_{j}\right)\leq\alpha_{\epsilon}$.

Conversely, if we take for initial measure $\nu$ such that its projection on the first layer is $\lambda_{1}=\chapeau 1$ and define $\nu_{\epsilon}^{t}=\left(\maj_{3}^{I,E}\right)_{\epsilon}^{t}\nu$. Denote by $\left(Y^{t},Z^{t}\right)$ a trajectory of $\nu_{\epsilon}^{t}$. We can bound the number of 1 in its trajectory by $\left(\maj_{3}^{I,E}\right)_{\epsilon}$ by another perturbation of $\maj_{3}^{I,E}$: define $\epsilon^{\prime}=\epsilon+\frac{\frac{1}{3}-\epsilon}{2}<\frac{1}{3}$. To conclude, its suffices to prove the following claim:
\begin{claim}
\label{claim:minoration_alpha_t}There exists $j\in\N$ such that if $\alpha^{t}\coloneqq P\left(Y_{j}^{t}=1\right)$, then $\alpha^{t}\geq h_{\epsilon^{\prime}}\left(\alpha^{t-1}\right)$ and $\alpha^{t}\geq1-\alpha_{\epsilon^{\prime}}$.
\end{claim}

If the claim is verified, then $\nu_{\epsilon}^{t}\left(\left[\left(1,0\right)\right]_{j}\cup\left[\left(1,1\right)\right]_{j}\right)\geq1-\alpha_{\epsilon^{\prime}}$. One can then observe that for $t\in\N$, 
\[
\mu_{\epsilon}^{t}\left(\left[\left(1,0\right)\right]_{j}\cup\left[\left(1,1\right)\right]_{j}\right)\leq\alpha_{\epsilon}<\frac{1}{2}<1-\alpha_{\epsilon^{\prime}}\leq\nu_{\epsilon}^{t}\left(\left[\left(1,0\right)\right]_{j}\cup\left[\left(1,1\right)\right]_{j}\right)
\]
and so $\mu_{\epsilon}^{t}$ and $\nu_{\epsilon}^{t}$ cannot converge towards the same measure, so $\left(\maj_{3}^{I,E}\right)_{\epsilon}$ is not ergodic. Moreover, their accumulation sets of the Ces\`aro mean of $\left(\mu_{\epsilon}^{t}\right)$ and $\left(\nu_{\epsilon}^{t}\right)$ are disjoint by the same argument, which proves the existence of at least two invariant measures, and so $\epsilon\notin\uni\left(\maj_{3}^{I,E}\right)$.
\begin{proof}[Proof of Claim \ref{claim:minoration_alpha_t}]
Let us first fix $j$. As $I_{3^{t}}=\left[\epsilon_{0}\pm\frac{1}{t}\right]$, one has $\epsilon\notin I_{3^{t}}$ for $t$ large enough. Similarly $\left|E_{3^{t}}\right|\tendsto t{+\infty}+\infty$ so by Bienaym\'e-Chebychev inequality, $P\left(2\cdot\text{mean}\left(Z^{t},i\right)\in I_{i}\right)\eqqcolon p_{i}\tendsto i{+\infty}0$. Let us fix $j$ large enough such that $p_{j}\frac{2g\left(1-\alpha_{\epsilon^{\prime}}\right)}{2g\left(1-\alpha_{\epsilon^{\prime}}\right)-1}\leq\frac{\frac{1}{2}\left(\frac{1}{3}-\epsilon\right)}{1-\epsilon}$ with $g:x\mapsto x^{3}+3x^{2}\left(1-x\right)$. Let us show that $\alpha^{t}\geq h_{\epsilon^{\prime}}\left(\alpha^{t-1}\right)$ and $\alpha^{t}\geq1-\alpha_{\epsilon^{\prime}}$.

By immediate induction on $k\in\left\llbracket 1,t-1\right\rrbracket $, $\left(Y_{i}^{k}\right)_{i\in\mathcal{N}_{j}^{t-k}}$ and $\left(Z_{E_{i}}^{k}\right)_{i\in\mathcal{N}_{j}^{t-k-1}}$ are families of i.i.d. random variables, independent from each other for each $k$ (the $Z_{i}^{k}$ are defined by the errors, which are independent, while for each $i\in\mathcal{N}_{j}^{t-k-1}$, $Y_{i}^{k+1}$ only depends on the error at $i$ at time $k+1$, on $\left(Y_{l}^{k}\right)_{l\in\mathcal{N}_{i}}$ and on $\left(Z_{l}^{k}\right)_{l\in E_{i}}$, which are independent from the others as $\mathcal{N}_{i}\subset\mathcal{N}_{j}^{t-k}$ and $E_{i}\subset E_{\mathcal{\mathcal{N}}_{j}^{t-k-1}}$ and those objects are disjoints by hypotheses). In particular, for $k=t-1$, $Y_{3j+1}^{t-1}$, $Y_{3j+2}^{t-1}$ and $Y_{3j+3}^{t-1}$ are i.i.d. and independent from $Z_{E_{j}}$. One then has 
\begin{eqnarray*}
\left\{ Y_{j}^{t}=1\right\} &=&\left\{ \text{Error at }j\text{ which creates a 1}\right\}\\
&& \cup\left(\left\{ 2\cdot\text{mean}\left(z,j\right)\notin I_{j}\right\} \cap\left\{ \text{maj}\left(Y_{3j+1}^{t-1},Y_{3j+2}^{t-1},Y_{3j+3}^{t-1}\right)=1\right\} \right)
\end{eqnarray*}

which leads to
\[
\alpha^{t}=\frac{\epsilon}{2}+\left(1-\epsilon\right)\left(1-p_{j}\right)g\left(P\left(Y_{3j+1}^{t-1}=1\right)\right).
\]
As $i\mapsto p_{i}$ is non-increasing, we easily obtain $P\left(Y_{3j+1}^{t-1}=1\right)\geq P\left(Y_{j}^{t-1}=1\right)$: so as $g:x\mapsto x^{3}+3x^{2}\left(1-x\right)$ is increasing, 
\[
\alpha^{t}\geq\frac{\epsilon}{2}+\left(1-\epsilon\right)\left(1-p_{j}\right)g\left(\alpha^{t-1}\right)\eqqcolon\Psi\left(\epsilon,p_{j},\alpha^{t-1}\right).
\]
For $\alpha>\frac{1}{2}$, $g\left(\alpha\right)>\frac{1}{2}$ and $\Psi\left(\epsilon,p,\alpha\right)\geq h_{\epsilon^{\prime}}\left(\alpha\right)\Leftrightarrow\epsilon^{\prime}\geq\epsilon+p\left(1-\epsilon\right)\frac{2g\left(\alpha\right)}{2g\left(\alpha\right)-1}$. We can then finish by induction on $t$: if $\alpha^{t-1}\geq1-\alpha_{\epsilon^{\prime}}$ (which is true at $t=1$), then $g\left(\alpha^{t-1}\right)>\frac{1}{2}$ and $\frac{2g\left(\alpha\right)}{2g\left(\alpha\right)-1}\leq\frac{2g\left(1-\alpha_{\epsilon^{\prime}}\right)}{2g\left(1-\alpha_{\epsilon^{\prime}}\right)-1}$, so by our choice of $p_{j}$ we get $\alpha^{t}\geq\Psi\left(\epsilon,p_{j},\alpha^{t-1}\right)\geq h_{\epsilon^{\prime}}\left(\alpha^{t-1}\right)$. As $h_{\epsilon^{\prime}}$ leaves $\left[1-\alpha_{\epsilon^{\prime}},1\right]$ stable, we then have $\alpha^{t}\geq1-\alpha_{\epsilon^{\prime}}$, and the claim is proven.
\end{proof}
\end{proof}
With a different choice of subsets $I_{3^{t}}=\left[0,b-\frac{1}{t}\right]$ with $b\leq\frac{1}{3}$, one can obtain a function of $\mathcal{C}\left(\A^{\N}\right)$ whose perturbation is ergodic only for $\epsilon\in\left[0,b\right)$ (when $\epsilon<\frac{1}{3}$). Similarly, with a choice of $I_{3^{t}}=\left[a+\frac{1}{t},b-\frac{1}{t}\right]$ the ergodicity is obtained only when $\epsilon\in\left(a,b\right)$. From there, any open set $O=\cup_{i\in\N}\left(a_{i},b_{i}\right)$ of $\left[0,\frac{1}{3}\right)$ can be obtained with $I_{3^{t}}=\bigcup_{i\leq t}\left[a_{i}+\frac{1}{t},b_{i}-\frac{1}{t}\right]$.

\subsection{Majority vote on $2n+1$ cells}\label{Sec.majorityVote}

To obtain non-ergodicity for error rate larger than the $\frac{1}{3}$ of $\maj_{3}$, we can simply increase the number of cells the majority vote uses.
\begin{defn}
Let $n\in\N$. Denote by $\maj_{2n+1}:\left\{ 0,1\right\} ^{\N}\to\left\{ 0,1\right\} ^{\N}$ the continuous function defined by 
\[
\maj_{2n+1}\left(x\right)_{i}\coloneqq\maj\left(x_{\left(2n+1\right)i+1},\dots,x_{\left(2n+1\right)i+2n+1}\right).
\]
For $\epsilon\in\left[0,1\right]$, define the polynomials 
\begin{align*}
g_{n}\left(X\right) & =\sum_{k=n+1}^{2n+1}\binom{2n+1}{k}X^{k}\left(1-X\right)^{2n+1-k}\\
P_{n}^{\epsilon} & =\frac{\epsilon}{2}+\left(1-\epsilon\right)g_{n}\left(X\right)-X.
\end{align*}
\end{defn}

\begin{prop}
\label{prop:maj2n+1}Let $n\in\N$, $\epsilon\in\left[0,1\right]$ and $\alpha\in\left[0,1\right]$. One has the following equivalences:
\begin{align*}
\lambda_{\alpha}\text{ is }\maj_{2n+1,\epsilon}\text{-invariant} & \Longleftrightarrow\alpha-\left(1-\epsilon\right)g_{n}\left(\alpha\right)=\frac{\epsilon}{2}.\\
 & \Longleftrightarrow\alpha\text{ is a root of }P_{n}^{\epsilon}
\end{align*}
\end{prop}

\begin{proof}
The proof is identical as the one for Proposition \ref{prop:maj3} in the case of $\maj_{3}$. For the direct direction, observe that if $B_{n}$ is a random variable with law $\text{Bin}\left(2n+1,\alpha\right)$, then $g_{n}\left(\alpha\right)=P\left(B_{n}\geq n+1\right)$. Therefore $\maj_{2n+1,\epsilon}\lambda_{\alpha}\left(\left[1\right]_{0}\right)=\frac{\epsilon}{2}+\left(1-\epsilon\right)g_{n}\left(\alpha\right)$. For the reciprocal, the local neighborhoods are $\mathcal{N}_{i}=\left\llbracket \left(2n+1\right)i+1,\left(2n+1\right)i+2n+1\right\rrbracket $ which are pairwise disjoints, which implies that the image of a Bernoulli measure by $\maj_{2n+1,\epsilon}$ is still a Bernoulli measure.
\end{proof}
One can make the following observation on the roots of $P_{n}^{\epsilon}$.
\begin{lem}
Let $n\in\N^{*}$ and $\epsilon>0$ such that $P_{n}^{\epsilon}$ admits a root $\alpha_{n,\epsilon}<\frac{1}{2}$. Then for all $0\leq\epsilon^{\prime}<\epsilon$, $P_{n}^{\epsilon^{\prime}}$ admits a root $\alpha_{n,\epsilon^{\prime}}<\frac{1}{2}$.
\end{lem}

\begin{proof}
For $\epsilon^{\prime}=0$, the result is immediate as $P_{n}^{0}$ admits $0$ as a root. Suppose now that $0<\epsilon^{\prime}<\epsilon$. Then 
\[
P_{n}^{\epsilon^{\prime}}=P_{n}^{\epsilon}+\left(\epsilon-\epsilon^{\prime}\right)\underset{<0\text{ on }\left[0,\frac{1}{2}\right[}{\underbrace{\left(g_{n}-\frac{1}{2}\right)}}.
\]
In particular, if one evaluates at $\alpha_{n,\epsilon}$, the result is $P_{n}^{\epsilon^{\prime}}\left(\alpha_{n,\epsilon}\right)<0$. As $P_{n}^{\epsilon^{\prime}}\left(0\right)>0$, the intermediate values theorem gives a root strictly smaller than $\frac{1}{2}$.
\end{proof}
The following lemma shows that by increasing the number of cells the majority vote takes place on, we can obtain maps whose perturbation stays non-ergodic for an arbitrary error rate.
\begin{lem}
Let $\epsilon<1$. Then there exists $N\in\N$ such that for all $n\geq N,$ $\maj_{2n+1,\epsilon}$ admits several invariant measures, i.e. $\epsilon\notin\uni\left(\maj_{2n+1}\right)$.
\end{lem}

\begin{proof}
Let $\alpha\in\left[0,\frac{1}{2}\right[$. Fix $B_{n}$ a random variable with law $\text{Bin}\left(2n+1,\alpha\right)$, with $n$ large enough such that
\[
n-\left(2n+1\right)\alpha=n\left(1-2\alpha\right)-\alpha>0.
\]
Then by Bienaym\'e-Chebychev inequality,
\begin{align*}
g_{n}\left(\alpha\right) & =P\left(B_{n}>n\right)\\
 & \leq P\left(\left|X_{n}-\left(2n+1\right)\alpha\right|>n-\left(2n+1\right)\alpha\right)\\
 & \leq\frac{\left(2n+1\right)\alpha\left(1-\alpha\right)}{\left(n\left(1-2\alpha\right)-\alpha\right)^{2}}\\
g_{n}\left(\alpha\right) & \leq\frac{\left(2n+1\right)\alpha\left(1-\alpha\right)}{n^{2}\left(1-2\alpha\right)^{2}}\tendsto n{\infty}0.
\end{align*}
Define $h_{n}:\alpha\longmapsto\alpha-\left(1-\epsilon\right)g_{n}\left(\alpha\right)$. $h_{n}$ is continuous, with $h_{n}\left(0\right)=0$ and $h_{n}\left(\frac{1}{2}\right)=\frac{\epsilon}{2}$. For all $\alpha<\frac{1}{2}$, the previous limit gives $h_{n}\left(\alpha\right)\tendsto n{\infty}\alpha$. In particular, for $\beta_{\epsilon}\coloneqq\frac{1}{2}\left(\frac{\epsilon}{2}+\frac{1}{2}\right)\in\left]\frac{\epsilon}{2},\frac{1}{2}\right[$, there exists $N\in\N$ such that for all $n\geq N$, $h_{n}\left(\beta_{\epsilon}\right)\in\left]\frac{\epsilon}{2},\frac{1}{2}\right[$. By intermediate values theorem, there exists $\alpha_{n,\epsilon}<\frac{1}{2}$ such that $h_{n}\left(\alpha_{n,\epsilon}\right)=\frac{\epsilon}{2}$.

Thus, for all $n\geq N$ $P_{n}^{\epsilon}$ admits at least two distinct roots: $\alpha_{\epsilon}$ and $\frac{1}{2}$. The Bernoulli measures associated to those parameters are distinct invariant measures, and so $\epsilon\notin\uni\left(\maj_{2n+1}\right)$.
\end{proof}
Proposition \ref{prop:maj2n+1} and the previous lemmas then lead to the following Corollary.
\begin{cor}
\label{cor:lnto1}For $n\in\text{\ensuremath{\N\ }}$, define $l_{n}\coloneqq\inf\left\{ \epsilon<1\mid\epsilon\in\uni\left(\maj_{2n+1}\right)\right\} $. Then $l_{n}\tendsto n{\infty}1$.
\end{cor}

\begin{rem}
For instance, we showed in Proposition \ref{prop:maj3} that $l_{1}\geq\frac{1}{3}$.
\end{rem}

As in the case of $\maj_{3}$, we can add an auxiliary layer where the maps acts as the constant equals to $0$. Then density of $1$ on a fixed sample gives a approximation of $\frac{\epsilon}{2}$ where $\epsilon$ is the error rate of the computer perturbation of the constructed map.
\begin{defn}
Let $I=\left(I_{i}\right)_{i\in\N}$ be a sequence of subsets of $\left[0,1\right]$, and $E=\left(E_{i}\right)_{i\in\N}$ a sequence of finite subsets of $\N$. For $x=\left(y,z\right)\in\left(\left\{ 0,1\right\} \times\left\{ 0,1\right\} \right)^{\N}$, define $x_{i}\coloneqq\left(y_{i},z_{i}\right)$. $\maj_{2n+1}^{I,E}\in \mathcal{C}\left(\A^{\N}\right)$ denotes the function with local neighborhood $\left\llbracket \left(2n+1\right)i+1,\left(2n+1\right)i+2n+1\right\rrbracket \cup E_{i}$ and local rule
\[
\maj_{2n+1}^{I,E}\left(x\right)_{i}\coloneqq\begin{cases}
\left(0,0\right) & \text{if }2\cdot\text{mean}\left(z,i\right)\in I_{i}\\
\left(\maj\left(y_{\left\llbracket \left(2n+1\right)i+1,\left(2n+1\right)i+2n+1\right\rrbracket }\right),0\right) & \text{otherwise}
\end{cases}
\]
with 
\[
\text{mean}\left(z,i\right)\coloneqq\frac{1}{\left|E_{i}\right|}\sum_{j\in E_{i}}z_{j}
\]
the frequency of symbols $1$ on the second layer of the cells in $E_{i}$.
\end{defn}

We can then show that with a good choice of $I$ and $E$, we can realize any open set of $\left[0,l_{n}\right[$ as an ergodicity set. We can also force the ergodicity of the perturbation for $\epsilon\geq l_{n}$, which leads to the following Proposition.
\begin{prop}
\label{prop:O=00005Ccup=00005Bln,1=00005D}Let $n\in\N^{*}$ and $O$ an open subset of $\left[0,1\right]$. There exists families of subsets $I$ and $E$ \textup{such that
\[
\uni\left(\maj_{2n+1}^{I,E}\right)=\erg\left(\maj_{2n+1}^{I,E}\right)=O\cup\left[l_{n},1\right].
\]
}
\end{prop}

\begin{proof}
For easier notations, define $N=2n+1$. One can decompose
\[
O\cup\left[l_{n},1\right]=\left(O\cap\left[0,l_{n}\right[\right)\cup\left[l_{n},1\right]
\]
and there exist sequences $a,b$ of $\left[0,l_{n}\right]$ such that $O\cap\left[0,l_{n}\right[=\bigcup_{i\in\N}\left]a_{i},b_{i}\right[$. The local neighborhood of $\left\{ 0\right\} $ by the continuous function $\maj_{N}^{t}$ is $\mathcal{N}^{t}=\left\llbracket \frac{N^{t}-1}{N-1},N^{t}+\frac{N^{t}-N}{N-1}\right\rrbracket \ni N^{t}$. As in the proof of Proposition \ref{prop:maj3IE}, we can choose $I$ and $E$ such that:
\begin{itemize}
\item $i\mapsto I_{i}$ is constant on each $N^{t}$, here with value $I_{N^{t}}=\left(\bigcup_{i\leq t}\left[a_{i}+\frac{1}{t},b_{i}-\frac{1}{t}\right]\right)\cup\left[l_{n}-\frac{1}{t},1\right]$.
\item $i\mapsto\left|E_{i}\right|$ is constant on each $N^{t}$, here with value $\left|E_{N^{t}}\right|=\left(N+1\right)^{t}$.
\item $E_{i}$ are pairwise disjoints.
\item For all $i\in\N$ and $t\geq1$, $\mathcal{N}_{i}^{t-1}\cap E_{\mathcal{N}_{i}^{t-2}}=\emptyset$.
\end{itemize}
We can define the $E_{i}$ by induction, beginning with $E_{0}=\left\{ N+1\right\} $. If $E_{i}=\left\llbracket a,b\right\rrbracket $ with $i\in\mathcal{N}^{k}$, define $E_{i+1}$ as:
\begin{itemize}
\item if $i+1\in\mathcal{N}^{k}$ as well, then $E_{i+1}=\left\llbracket b+1,b+\left(N+1\right)^{k}\right\rrbracket $.
\item if $i+1\in\mathcal{N}^{k+1}$, then $E_{i+1}=\left\llbracket b^{\prime}+1,b^{\prime}+\left(N+1\right)^{k+1}\right\rrbracket $ with $b^{\prime}\coloneqq\max\left(b,\max\left(\mathcal{N}^{k+2}\right)\right)$. 
\end{itemize}
By construction, the $E_{i}$ verify the previously listed properties. 

The end of the proof is analogous to the one of Proposition \ref{prop:maj3IE}, using for a fixed $\epsilon$ 
\[
p_{i}\coloneqq P\left(2\cdot\text{mean}\left(Z^{1},i\right)\in I_{i}\right)
\]
with $Z^{t}$ the (common) second layer of a coupling $\left(X^{t},Y^{t}\right)$ of two initial measures $\mu$ and $\nu$. 
\begin{itemize}
\item For $\epsilon\in O\cup\left[l_{n},1\right]$, one can verify that $\left[\epsilon\pm t\right]\subset I_{N^{t}}$ for $t$ large enough, and thus the Bienaym\'e-Chebychev inequality leads to
\[
p_{N^{t}}\geq1-\frac{t^{2}}{\left|E_{N^{t}}\right|}=1-\frac{t^{2}}{\left(N+1\right)^{t}},
\]
where $P\left(X_{0}^{t}=Y_{0}^{t}\right)\geq\left(p_{N^{t-1}}\right)^{N^{t-2}}\tendsto t{\infty}1$. One can likewise show that $P\left(X_{j}^{t}=Y_{j}^{t}\right)\tendsto t{\infty}1$ for all $j\in\N$ and the ergodicity of $\left(\maj_{N}^{I,E}\right)_{\epsilon}$ follows.
\item For $\epsilon<l_{n}$ and $\epsilon\notin O$, one has $\left[\epsilon\pm t\right]\cap I_{N^{t}}=\emptyset$ for $t$ large enough and therefore $p_{i}\tendsto i{\infty}0$. The same analysis show that using a initial measure $\mu$ whose projection on the second layer is $\lambda_{0}$, one has 
\[
\left(\maj_{N}^{I,E}\right)_{\epsilon}^{t}\mu\left(\left[1,0\right]_{j}\cup\left[1,1\right]_{j}\right)\leq\alpha_{n,\epsilon}
\]
where $\alpha_{n,\epsilon}$ is the root strictly smaller than $\frac{1}{2}$ of $P_{n}^{\epsilon}$, observing that $h_{n,\epsilon}\left(\alpha\right)\coloneqq\frac{\epsilon}{2}+\left(1-\epsilon\right)g_{n}\left(\alpha\right)$ leaves $\left[0,\alpha_{N,\epsilon}\right]$ and $\left[1-\alpha_{n,\epsilon},1\right]$ stable. Conversely, using $\nu$ a initial measure whose projection on the second coordinate is $\lambda_{1}$ leads, for $j$ large enough, to
\[
1-\alpha_{n,\epsilon^{\prime}}\leq\left(\maj_{N}^{I,E}\right)_{\epsilon}^{t}\nu\left(\left[1,0\right]_{j}\cup\left[1,1\right]_{j}\right)
\]
where $\epsilon^{\prime}=\frac{\epsilon+l_{n}}{2}<l_{n}$, and $\alpha_{n,\epsilon^{\prime}}$ the root strictly smaller than $\frac{1}{2}$ of $P_{n}^{\epsilon^{\prime}}$. The trajectories of $\mu$ and $\nu$ cannot have an accumulation point in common, and thus $P_{n}^{\epsilon}$ admits at least two different invariant measures. 
\end{itemize}
\end{proof}

\subsection{Construction of $M$ and $M^{I,E}$}\label{Sec.RealisationOpen}

We succeeded in controlling the ergodicity of $\maj_{2n+1}$ in the case $\epsilon<l_{n}$ with $l_{n}$ converging to $1$, and forcing the ergodicity on $\left[l_{n},1\right]$. The last step consist in simulating all those maps independently on different lines of $\N^{2}$, and then use a continuous bijection of $\N^{2}$ to $\N$ to transform the continuous function of $\A^{\N^{2}}$ into a continuous function of $\A^{\N}$.
\begin{defn}
Denote by $\widetilde{M}\in C\left(\left\{ 0,1\right\} ^{\N^{2}}\right)$ the function applying the dynamics of $\maj_{2j+1}$ on the $j^{\text{th}}$ row: for $x\in\left\{ 0,1\right\} ^{\N^{2}}$, 
\[
\widetilde{M}\left(x\right)_{i,j}=\left(\maj_{2j+1}\left(x_{\cdot,j}\right)\right)_{i}.
\]
\end{defn}

Define $\varphi:\N^{2}\to\N$ by $\varphi\left(i,j\right)=2^{j}-1+i2^{j+1}$: it sends the $j^{\text{th}}$ row of $\N^{2}$ on $2^{j}-1+2^{j+1}\N$. It defines the (continuous) change of coordinates $\Phi:\A^{\N}\to\A^{\N^{2}}$ via $\Phi\left(x\right)_{i,j}\coloneqq x_{\varphi\left(i,j\right)}$.
\begin{defn}
Denote by $M\in C\left(\left\{ 0,1\right\} ^{\N}\right)$ the function defined by
\[
M\coloneqq\Phi^{-1}\circ\widetilde{M}\circ\Phi.
\]
\end{defn}

It corresponds to sending the dynamics $\maj_{2j+1}$ on $\A^{\left(2^{j}-1\right)+2^{j+1}\N}$ ($\maj_{1}$ acts on the even cells, $\maj_{3}$ on the cells of $1+4\N$, etc.). 
\begin{figure}
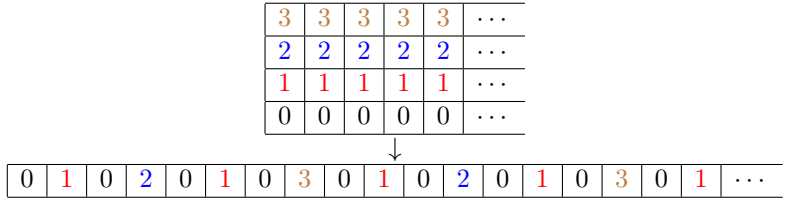

\begin{centering}
\begin{tabular}{|c|c|c|c|c|c}
\hline 
\textcolor{brown}{3} & \textcolor{brown}{3} & \textcolor{brown}{3} & \textcolor{brown}{3} & \textcolor{brown}{3} & $\cdots$\tabularnewline
\hline 
\textcolor{blue}{2} & \textcolor{blue}{2} & \textcolor{blue}{2} & \textcolor{blue}{2} & \textcolor{blue}{2} & $\cdots$\tabularnewline
\hline 
\textcolor{red}{1} & \textcolor{red}{1} & \textcolor{red}{1} & \textcolor{red}{1} & \textcolor{red}{1} & $\cdots$\tabularnewline
\hline 
0 & 0 & 0 & 0 & 0 & $\cdots$\tabularnewline
\hline 
\end{tabular}
\par\end{centering}
\begin{centering}
$\downarrow$
\par\end{centering}
\begin{centering}
\begin{tabular}{|c|c|c|c|c|c|c|c|c|c|c|c|c|c|c|c|c|c|c}
\hline 
0 & \textcolor{red}{1} & 0 & \textcolor{blue}{2} & 0 & \textcolor{red}{1} & 0 & \textcolor{brown}{3} & 0 & \textcolor{red}{1} & 0 & \textcolor{blue}{2} & 0 & \textcolor{red}{1} & 0 & \textcolor{brown}{3} & 0 & \textcolor{red}{1} & $\cdots$\tabularnewline
\hline 
\end{tabular}
\par\end{centering}
\centering{}\caption{Illustration of the partition of $\protect\N$ from $\protect\N^{2}$. $M$ acts as $\protect\maj_{1}$ on the cells denoted by 0, as $\protect\maj_{3}$ on the cells denoted by \textcolor{red}{1}, as $\protect\maj_{5}$ on the cells denoted by \textcolor{blue}{2} and as $\protect\maj_{7}$ on the cells denoted by \textcolor{brown}{3}.}
\end{figure}

\begin{prop}
\[
\uni\left(M\right)=\erg\left(M\right)=\left\{ 1\right\} .
\]
\end{prop}

\begin{proof}
The respective computer perturbation of $\widetilde{M}$ and $M$ are also conjugated by $\Phi$:
\[
M_{\epsilon}=\Phi^{-1}\circ\widetilde{M}_{\epsilon}\circ\Phi
\]
in the sense of the probability kernels: for all $x\in\left\{ 0,1\right\} ^{\N}$ and $A\in\mathfrak{B}\left(\left\{ 0,1\right\} ^{\N}\right)$, 
\[
M_{\epsilon}\left(x,A\right)=\widetilde{M}_{\epsilon}\left(\Phi\left(x\right),\Phi\left(A\right)\right).
\]
As $\Phi:\A^{\N}\to\A^{\N^{2}}$ is a bijection, one easily has $\uni\left(M\right)=\uni\left(\widetilde{M}\right)$ and $\erg\left(M\right)=\erg\left(\widetilde{M}\right)$. Therefore, we can show the result on the ergodicity sets of $\widetilde{M}$.

The computer perturbation $\widetilde{M}_{\epsilon}$ of $\widetilde{M}$ acts as the computer perturbation of each $\maj_{2j+1}$, independently: $\widetilde{M}_{\epsilon}$ acts as $\maj_{2j+1,\epsilon}$ on the row $j$. In particular, if $\pi_{\epsilon}^{j}$ is a $\maj_{2j+1,\epsilon}$-invariant measure, the product measure $\bigotimes_{j\in\N}\pi_{\epsilon}^{j}$ is a $\widetilde{M}_{\epsilon}$-invariant measure, which implies $\erg\left(\widetilde{M}\right)\subset\uni\left(\widetilde{M}\right)\subset\bigcap_{j\in\N}\uni\left(\maj_{2j+1}\right)$.

By Corollary \ref{cor:lnto1}, for all $\epsilon<1$ there exists a row $j$ for $j$ large enough such that $\maj_{2j+1,\epsilon}$ admits several invariant measures. Therefore for all $\epsilon<1$, $\widetilde{M}_{\epsilon}$ admits several invariant measures, and $\uni\left(\widetilde{M}\right)\subset\left\{ 1\right\} $. At $\epsilon=1$ the ergodicity is guaranteed by the independence on the noise (true for all computer perturbation). Thus $\left\{ 1\right\} \subset\erg\left(\widetilde{M}\right)\subset\uni\left(\widetilde{M}\right)\subset\left\{ 1\right\} $ and the result is proven.
\end{proof}
We have then a fairly simple continuous map which admits several invariant measures for any error rate (except $1$). To force its ergodicity on any open subset $O$, we add similarly as before a auxiliary layer where the constant function equals to $0$ is simulated, to have at disposal an approximation of the error rate. This time, the families $I$ and $E$ are indexed by $\N^{2}$, as we can take different families for each line of $\N^{2}$ to re-use the previous results .
\begin{defn}
Let $I=\left(I_{i,j}\right)$ and $E=\left(E_{i,j}\right)$ be families of subsets of $\left[0,1\right]$ and $\N$. For $x=\left(y,z\right)\in\left(\left\{ 0,1\right\} \times\left\{ 0,1\right\} \right)^{\N^{2}}$, define $x_{i}=\left(y_{i},z_{i}\right)$. Denote by $\widetilde{M}^{I,E}$ the function defined by
\[
\widetilde{M}^{I,E}\left(x\right)_{i,j}\coloneqq\left(\maj_{2j+1}^{I_{j},E_{j}}\left(x_{\cdot,j}\right)\right)_{i}
\]
with $I_{j}=\left(I_{i,j}\right)_{i\in\N}$ and $E_{j}=\left(E_{i,j}\right)_{i\in\N}$. We make it a continuous function of $\left(\left\{ 0,1\right\} ^{2}\right)^{\N}$ by $M^{I,E}=\Phi^{-1}\circ\widetilde{M}^{I,E}\circ\Phi$.
\end{defn}

\begin{prop}
\label{prop:O=00005Ccup=00007B1=00007D}Let $O$ be a open subset of $\left[0,1\right]$. For all $j\in\N$, define $I_{j}$ and $E_{j}$ the families corresponding to the ones defined in Proposition \ref{prop:O=00005Ccup=00005Bln,1=00005D} to obtain $O\cup\left[l_{j},1\right]$ as ergodicity sets of $\maj_{2j+1}^{I_{j},E_{j}}$. Then 
\[
\erg\left(M^{I,E}\right)=\uni\left(M^{I,E}\right)=\bigcap_{j\in\N}\uni\left(\maj_{2j+1}^{I_{j},E_{j}}\right)=\bigcap_{j\in\N}\erg\left(\maj_{2j+1}^{I_{j},E_{j}}\right)=O\cup\left\{ 1\right\} .
\]
\end{prop}

\begin{proof}
By the choice of $I$ and $E$, one has
\begin{align*}
\bigcap_{j\in\N}\uni\left(\maj_{2j+1}^{I_{j},E_{j}}\right) & =\bigcap_{j\in\N}\erg\left(\maj_{2j+1}^{I_{j},E_{j}}\right)\\
 & =\bigcap_{j\in\N}O\cup\left[l_{j},1\right]\\
 & =O\cup\left\{ 1\right\} 
\end{align*}
since $l_{j}\tendsto j{\infty}1$ by Corollary \ref{cor:lnto1}. As in the precedent proof, one easily show that $\erg\left(\widetilde{M}^{I,E}\right)=\erg\left(M^{I,E}\right)$, and the same for $\uni$. The same argument of product measures leads to
\[
\erg\left(\widetilde{M}^{I,E}\right)\subset\uni\left(\widetilde{M}^{I,E}\right)\subset O\cup\left\{ 1\right\} .
\]

Conversely, let $\epsilon\in O\cup\left\{ 1\right\} =\bigcap_{j\in\N}\erg\left(\maj_{2j+1}^{I_{j},E_{j}}\right)$. We can once again use a coupling $\left(X^{t},Y^{t}\right)$ of the trajectories of two initial measures $\mu$ and $\nu$ of $\left(\left\{ 0,1\right\} \times\left\{ 0,1\right\} \right)^{\N^{2}}$ by $\left(\widetilde{M}^{I,E}\right)_{\epsilon}$. For all $j\in\N$, by independence of the perturbations, the projections on the $j^{\text{th}}$ row $P_{j}\left(X^{t}\right)$ and $P_{j}\left(Y^{t}\right)$ are a coupling of the trajectories$P_{j}\mu$ and $P_{j}\nu$ by $\left(\maj_{2j+1}^{I_{j},E_{j}}\right)_{\epsilon}$. As $\epsilon\in\erg\left(\maj_{2j+1}^{I_{j},E_{j}}\right)$, we can conclude that for all $i,j\in\N$
\[
P\left(X_{i,j}^{t}=Y_{i,j}^{t}\right)=P\left(\left(P_{j}\left(X^{t}\right)\right)_{i}=\left(P_{j}\left(Y^{t}\right)\right)_{i}\right)\tendsto t{\infty}1.
\]
Therefore $O\cup\left\{ 1\right\} \subset\erg\left(\widetilde{M}^{I,E}\right)$ and the result is proven.
\end{proof}

\subsection{Realizing any $G_{\delta}$ set}\label{Sec.ConstructionGdelta}

The method consisting in transforming a continuous function of $\A^{\N^{2}}$ into a continuous function of $\A^{\N}$ can be used to obtain any open subset as an ergodicity set, by describing it as a countable intersection of sets we already could get ($O\cup\left\{ 1\right\} =\bigcap_{j}O\cup\left[l_{j},1\right]$). As a $G_{\delta}$ set is a countable intersection of open sets, the same method can be applied to obtain them.
\begin{thm}
\label{thm:Ergodicit=0000E9GDelta} Let $G$ be a $G_{\delta}$ subset of $\left[0,1\right]$ with $1\in G$. There exists $F$ a continuous function on $\left(\left\{ 0,1\right\} \times\left\{ 0,1\right\} \right)^{\N}$ such that the computer perturbation $F_{\epsilon}$ is ergodic if and only if $\epsilon\in G$. Otherwise, $F_{\epsilon}$ admits several invariant measures. Formally, 
\[
\uni\left(F\right)=\erg\left(F\right)=G.
\]
\end{thm}

\begin{proof}
One can describe $G$ as $G=\bigcap_{j}\left(O_{j}\cup\left\{ 1\right\} \right)$, with $\left(O_{j}\right)$ a sequence of open subsets of $\left[0,1\right]$. For each $n\in\N$, choose $M_{j}\in C\left(\left(\left\{ 0,1\right\} \times\left\{ 0,1\right\} \right)^{\N}\right)$ a continuous function verifying $\erg\left(M_{j}\right)=\uni\left(M_{j}\right)=O_{j}\cup\left\{ 1\right\} $ (for instance, the function $M^{I,E}$ described in Proposition \ref{prop:O=00005Ccup=00007B1=00007D}). Define then $\widetilde{F}$ the continuous function of $\left(\left\{ 0,1\right\} \times\left\{ 0,1\right\} \right)^{\N^{2}}$ simulating the dynamics of $M_{j}$ on the $j^{\text{th}}$ row:
\[
\widetilde{F}\left(x\right)_{i,j}\coloneqq\left(M_{j}\left(x_{\cdot,j}\right)\right)_{i}.
\]
The final function $F$ is then $F\coloneqq\Phi^{-1}\circ\widetilde{F}\circ\Phi$. By an identical proof as for Proposition \ref{prop:O=00005Ccup=00007B1=00007D} (coupling for ergodicity, and product measures for the existence of several invariant measures), one gets

\[
\erg\left(F\right)=\uni\left(F\right)=\bigcap_{j\in\N}\erg\left(M_{j}\right)=\bigcap_{j\in\N}\left(O_{j}\cup\left\{ 1\right\} \right)=G.
\]
\end{proof}

\section{Generic phase diagram}\label{sec.Generic}

We are interested in typical behavior. First, we consider the distribution of functions that realize a given $G_\delta$ set. By modifying the construction carried out in the previous section, we can demonstrate that the set of continuous functions realizing a given $G_\delta$ set containing $1$ but not $0$ is dense.

\begin{prop}\label{prop.densite}
Let $G$ be a $G_{\delta}$ subset of $\left[0,1\right]$ with $1\in G$ and $0\notin G$. The following set is dense in $(\mathcal{C}\left(\A^{\N}\right),d_{\infty})$:
$$\left\{F\in\mathcal{C}\left(\A^{\N}\right):\uni\left(F\right)=G\right\}.$$
\end{prop}
\begin{proof}
Let $F\in\mathcal{C}\left(\A^{\N}\right)$ and consider the function $F_G$ obtained by Theorem~\ref{thm:Ergodicit=0000E9GDelta} where $\uni\left(F^G\right)=G$ . Given $n\in\N$, since $F$ is continuous on the compact set $\A^\N$, there exists $r\geq n$ such that $F(x)_{[0,n]}$ depends only on $x_{[0,r]}$ for any $x\in\A^\N$. We construct the function $F'$, defined for $x\in\A^\N$, as follows:
\begin{eqnarray*}
 F'(x)_{[0,n]}&=&F(x)_{[0,n]}\\
 F'(x)_{[n+1,r]}&=&0^{r-n}\\
 F'(x)_{r+1}&=&x_{r+1}\\
 F'(x)_{[r+2,\infty]}&=&F^G(x_{[r+2,\infty]})
\end{eqnarray*}

By construction we deduce that $d_{\infty}(F,F')\leq 2^{-n}$. We then observe that the coordinates $[0,r]$, $\{r+1\}$ and $[r+2,\infty]$ are independent. On  $[0,r+1]$, the function $F_\epsilon$ acts as a Markov chain with a finite number of states and a positive rate for $\epsilon>0$ so it admits only one invariant measure. On $[r+2,\infty]$ the number of invariant measures is the same as that of $F^G_\epsilon$. We therefore conclude that $$\uni\left(F\right)\setminus\{0\}=\uni\left(F^G\right)\setminus\{0\}=G\setminus\{0\}=G.$$

For $\epsilon=0$, the coordinate $r+1$ produces at least two invariant measures. Thus $\uni\left(F\right)=G$.
\end{proof}

If $0\in G$ then $\left\{F\in\mathcal{C}\left(\A^{\N}\right):\uni\left(F\right)=G\right\}$ is not dense. Indeed, if $F(x)_i=x_i$, i.e. if $F$ fixes a coordinate, then there exists a neighborhood of $F$ in the space of continuous functions whose elements also fix this coordinate. Thus, the elements in this neighborhood have several invariant measures. Therefore, this neighborhood does not intersect the set of functions in question. 

We would like to characterize $\uni\left(F\right)$ for a generic continuous function $F$. A set is \emph{generic}, or \emph{co-meager}, if it contains a dense $G_\delta$ set. We have the following result.

\begin{thm}
The following set is generic in $(\mathcal{C}\left(\A^{\N}\right),d_{\infty})$:
$$\left\{F\in\mathcal{C}\left(\A^{\N}\right):\uni\left(F\right)=]0,1]\right\}.$$
\end{thm}
\begin{proof}
Denote $A_0=\left\{F\in\mathcal{C}\left(\A^{\N}\right):\uni\left(F\right)\setminus\{0\}=]0,1]\right\}$. We have $A_0=\bigcap_kA_k$ where
$$A_k=\left\{F\in\mathcal{C}\left(\A^{\N}\right): \forall\epsilon\in\left[\frac{1}{k},1\right],\ \text{diam}\left(\M^F_{\epsilon}\right)<\frac{1}{k}\right\}.$$

It is sufficient to prove that $A_k$ is an open set, the density of $A_0$ is obtained by Proposition~\ref{prop.densite}. Let $(F^i)_{i\in\N}$ be a sequence of element of $A_k^c$ which converges to $F$. For each $i$ there exists an element of $\left[\frac{1}{k},1\right]$ denoted $\epsilon_i$ such that the diameter of the set $M^{F^i}_{\epsilon_i}$ is greater than or equal to $\frac{1}{k}$. By compactness, we can assume that the sequence $(\epsilon_i)_{i\in\N}$ converges to $\epsilon$. Using the continuity Lemma~\ref{lem:ContinuityLemma}, we deduce that $\text{diam}\left(\M^{F}_{\epsilon}\right)\geq\frac{1}{k}$. So $F\in A_k^c$ and thus $A_k$ is open. 

Consider $\mathcal{U}$ the set of continuous function which fix at least a coordinate, that is to say
$$\mathcal{U}=\left\{F\in\mathcal{C}\left(\A^{\N}\right): \exists i\in\N,\ \forall x\in\A^\N,\  F(x)_i=x_i \right\}.$$
The set $\mathcal{U}$ is a dense open set such that $0\notin\uni\left(F\right)$ for any $F\in\mathcal{U}$. To conclude we just verify that the set considered in the theorem contains the intersection of $A_0$ with $\mathcal{U}$.
\end{proof}

\section{Some additional obstructions for the class of cellular automata}

Returning to the class of cellular automata, it is still open to characterize the possible sets written as $\uni\left(F\right)$.  In this case, other restrictions must be taken into consideration. Firstly, for a cellular automaton $F$, there exists $r_F$ such that $[r_F,1]\subset\uni\left(F\right)$ (see~\cite{MST19}). 

Secondly, the countable number of cellular automata leads to combinatorial constraints on the sets that can be realised as $\uni\left(F\right)$, in addition to the topological constraints of Proposition~\ref{prop:uni_est_un_G_delta}.  Since cellular automata can be considered a model of computation, we search for computable obstructions, as discussed in~\cite{HS18,GST23,Marsan-Sablik-2024}. It is possible to prove that if $F$ is computable then $\uni\left(F\right)$ is a $\Pi_2$-computable set in the sense that there exists a computable map $h:\Q^{2}\times\N^{2}\to\left\{ 0,1\right\} $ verifying for all $a\leq b\in\left[0,1\right]\cap\Q$,
\[
\left[a,b\right]\cap\uni\left(F\right)\neq\emptyset\Longleftrightarrow\forall k\in\N,\exists l\in\N,h\left(a,b,k,l\right)=1.
\]

We remark that the construction of Theorem \ref{thm:Ergodicit=0000E9GDelta} is itself computable, as long as the open sets $O_{j}$ used are uniformly $\Sigma_1$-computable, that is to say that we can use a Turing machine to obtain the inner approximation as $I_{N^{t}}=\bigcup_{i\leq t}\left[a_{i}+\frac{1}{t},b_{i}-\frac{1}{t}\right]$). Thus any $\Pi_2$-computable $G_\delta$ of $[0,1]$ which contains $1$ can be obtained as $\uni\left(F\right)$ where $F$ is a computable map.

\end{document}